\documentclass[10pt]{article}

\usepackage{amsfonts,amsmath,latexsym,amssymb}    % need for subequations

\title{{\smc fractional brownian flows}}

\author{{\smc Sreekar Vadlamani}\thanks{Research supported in part
by the  US-Israel Binational
Science Foundation, grant 2004064.}}

\newtheorem{lemma}{Lemma}[section]
\newtheorem{theorem}{Theorem}[section]
\newtheorem{remark}[theorem]{Remark}

\newcommand{\beq}{\begin{eqnarray}}
\newcommand{\eeq}{\end{eqnarray}}
\newcommand{\beqq}{\begin{eqnarray*}}
\newcommand{\eeqq}{\end{eqnarray*}}
\newcommand{\qed}{\hfill \ensuremath{\Box}}
\newcommand{\be}{\begin{equation}}
\newcommand{\ee}{\end{equation}}
\newcommand{\bea}{\begin{eqnarray}}
\newcommand{\eea}{\end{eqnarray}}
\newcommand{\beaa}{\begin{eqnarray*}}
\newcommand{\eeaa}{\end{eqnarray*}}
\newcommand{\beas}{\begin{eqnarray*}}
\newcommand{\eeas}{\end{eqnarray*}}
\newcommand{\ben}{\begin{enumerate}}
\newcommand{\een}{\end{enumerate}}
\newcommand{\bi}{\begin{itemize}}
\newcommand{\ei}{\end{itemize}}
\newcommand{\ba}{\begin{array}}

\newcommand{\ea}{\end{array}}

\newcommand{\eps}{\varepsilon}

\newcommand{\al}{\alpha}

\def\complex{\mathop{\raise .45ex\hbox{${\bf\scriptstyle{|}}$}
     \kern -0.40em {\rm \textstyle{C}}}\nolimits}

\newcommand{\smc}{\scshape}

\newcommand{\calH}{{\cal H}}

\newcommand{\iy}{\infty}

\begin{document}
\baselineskip=12pt

\maketitle

%------------------------BEGIN MAIN---------------------------------

%------------------------ABSTRACT-----------------------------------
\begin{abstract}

We consider stochastic flow on $\mathbb{R}^n$ driven by fractional
Brownian motion with Hurst parameter $H\in (\frac{1}{2},1)$, and
study tangent flow and the growth of the Hausdorff measure of  sub-manifolds
of $\mathbb{R}^n$ as they evolve under the flow.

The main result is a bound on the rate of (global) growth in terms of
the  (local) H\"older norm of the flow.

\end{abstract}

%-------------------------------------------------------------------

\section{Introduction}
\label{section.intro}

Our  main objective is to study the global geometric
properties of a manifold embedded in Euclidean space, as it evolves
under a stochastic flow of diffeomorphisms driven by a non diffusive
process. This follows on from our previous paper \cite{VaAd06} in which we
obtained precise estimates for the rate of growth of the {\em
Lipschitz-Killing curvatures}\footnote{For a detailed exposition on
Lipschitz-Killing curvatures, we refer the reader to
\cite{AdTay07}.} of  smooth, $(n-1)$-dimensional manifolds embedded
in $\mathbb{R}^n$, as they evolve under an isotropic Brownian flow.
%In the literature, Lipschitz-Killing curvatures are also referred to
%as {\em intrinsic volumes} or {\em curvature measures}, for they are
%invariant of the ambient space in which the manifold is embedded,
%and they are also invariant under rigid motions.
In this paper, however, we turn to the non-Markovian, non-diffusive
situation in which the flow is driven by fractional Brownian motions.

Although extensive literature is available for stochastic flows
driven by standard Brownian motion (see \cite{BaxHar86,Kun90}),
very little is known when the driver of the flow is changed to a
non-Markovian, non-diffusive process, such as fractional Brownian motion.
For instance, some of the very basic results concerning the tangent flow
are yet to be unearthed in the case when the flow is driven by fractional
Brownian motion. Here we intend to target precisely this aspect of the flow
on our way to the main result of this paper.

%As a natural extension of \cite{VaAd06}, one would want to study
%similar aspects of the flow when it is driven by a non-markovian
%process, resulting in a non-diffusive flow.
%obtain similar results on
%the rate of growth of some global geometric characteristics of a
%smooth manifold embedded in Euclidean space, as it evolves under a
%non-diffusive flow.

Recall that a fractional Brownian motion $\{B^H(t),\; t\ge 0\}$ with Hurst
parameter $H\in (0,1)$, is the zero mean Gaussian process with
covariance function
\be \label{equation.fBM.cov} E[B^H(s)B^H(t)] =
\frac{1}{2}\big(t^{2H} +s^{2H} -|t-s|^{2H}\big). \ee
When $H=1/2$, $B^H$ is the standard Brownian motion, which
is a Markov process and also a martingale. However for $H\neq 1/2$,
$B^H$ is neither a Markov process, nor a semi-martingale.

%rendering a need to develop stochastic calculus for processes which
%are not semi-martingales. This, and more, forms the core of
%\cite{AloMazNua01,CarCouMon03,DaiHey96,DecUst99,Lin95,Lyon98,Zah98},
%where various ways of defining integrals with respect to fractional
%Brownian motion are discussed and the corresponding stochastic
%calculus is developed.

%The fact that the driving process possesses stationary increments
%implies that the corresponding {\em noise} process is a stationary
%process, and so one can hope to obtain ergodic estimates.

%The primary questions related to the study of the various ways of
%defining stochastic integrals with respect to Gaussian processes, in
%particular fractional Brownian motions with $H> 1/2$ have been
%studied in detail by various researchers in .................. Two
%of the most natural ways being {\em pathwise definition} and using
%{\em Malliavin calculus}.

In order to construct a non-diffusive flows, we start with a
collection of independent fractional Brownian motions,
$\{B^H_{\gamma}\}_{\gamma\in \mathbb{N}}$, a collection  $\{U_{\gamma}\}_{\gamma\in
\mathbb{N}}$  of deterministic vector fields
on $\mathbb{R}^n$, and define, for some fixed
but generic set $I\subset\mathbb{N}$ with  $|I|<\iy$, where
$|I|$ denotes the cardinality of  $I$,
\be\label{equation.fBM.vector.fields}U_I(x,t) = \sum_{\gamma\in I}
U_{\gamma}(x)B^H_{\gamma}(t).\ee
%and where %For the moment, we do not impose any conditions
%on the vector fields, although we shall add various conditions as
%and when required.

The flow of diffeomorphisms $\Phi_t:\mathbb{R}^n\to \mathbb{R}^n$,
or, equivalently, the
stochastic flow driven by a fractional Brownian motion,
can then be defined pointwise  by setting
\be \label{equation.fBM.flow.candidate} \Phi_t(x)=x+\sum_{\gamma\in
I}``\int_0^t U_{\gamma}(\Phi_s(x))\;dB^H_{\gamma}(s)". \ee

Clearly, we shall need to place conditions on the vector fields for the
result to give a diffeomorphism, but, prior to that, we need to make sense of
the stochastic integrals here.

For $H=\frac{1}{2}$, the integral can be interpreted in  either the
 It\^{o} or a Stratonovich sense.
When $H\neq\frac{1}{2}$ the standard semimartingale
arguments cease to work and  we have to make a choice of definition.
There is a plethora of literature
on various ways to define an integral $\int_a^b f(s)\;dB^H(s)$,
where $f$ is random and $B^H$ the fractional
Brownian motion. See, for instance,
\cite{AloMazNua01,CarCouMon03,DecUst99,DaiHey96,Lin95,Zah98}.

We shall adopt the pathwise definition given by Z\"{a}hle
\cite{Zah98,Zah99}, based on which
Nualart and R{\u{a}}{\c{s}}canu (\cite{NuaRas02})
proved existence and uniqueness of the solutions of multidimensional
stochastic differential equations of the form
$$X_t=X_0+\int_0^t\sigma(s,X_s)\;dB^H(s)+\int_0^tb(s,X_s)\;ds,$$
for  $H>\frac{1}{2}$.
Using this,  Decreusefond and Nualart in \cite{DecNua06} established the
existence of a homeomorphic stochastic flow driven by fractional
Brownian motion, and so our flows are well defined. We note here that
stochastic integrals can also be defined for $H<1/2$ using Malliavin
calculus (see \cite{DecUst99}), but existence of solution of stochastic
integral equations of above type is not ensured.

Note that in
\eqref{equation.fBM.flow.candidate} we do not have a drift part, as
we intend to study flows driven purely by noise, which goes hand in hand
with the way stochastic flows have been defined in \cite{Kun90}.

%Around the same time, Z\"{a}hle \cite{Zah98,Zah99}, defined pathwise
%integrals of the form $\int_0^t u_sdB^H(s)$, for $H>\frac{1}{2}$,
%where $u$ and $B^H$ were considered to be elements of fractional
%Sobolev space. Apart from proving an It\^{o} type formula,
%connections were established between the pathwise definition of the
%stochastic integral with various other interpretations of integrals.
%Elsewhere, Ruzmaikina \cite{Ruz00} independently, obtained similar
%results by approximating the stochastic integral by a
%Riemann-Stieltjes sum.

Now that it is clear, in principle, which flows we are considering, we can
turn to the geometry. Consider a fixed  $m$-dimensional, ($m<n$) $C^2$  manifold
embedded in $\mathbb{R}^n$, and consider its image under $\Phi$, setting
$$
M_t = \Phi_t(M) = \{x\in \mathbb{R}^n: x=\Phi_t(y)\ \text{for some}\ y\in M\}.
$$
Our interest is how $M_t$ behaves as a function of $t$.

 Although in  \cite{VaAd06} we were able to obtain information on all
the Lipschitz-Killing curvatures of $M_t$, in the current, non-diffusion
scenario everything is much harder, and so we shall suffice by studying only
the size of $M_t$, as measured through its $m$-dimensional Hausdorff measure,
$\calH_m(M_t)$, which basically measures the $m$-dimensional Lebesgue measure
of the set $M_t$. Our main result is Theorem \ref{theorem.fBM.growth.Hausdorff.measure},
however one can see the main flavour of the result already for a flow driven
by a single fractional Brownian motion. In this case we have

\begin{theorem}
\label{theorem.growth.Hausdorff.1}
In the notation above, assuming that $|I|=1$ in
\eqref{equation.fBM.flow.candidate}, and under conditions
(A1)--(A3) of Section 2 on the vector field $U$, for every $\beta <H$ and $H>\frac12$
%
%Let $M_t$ be the image under the
%fractional flow $\Phi_t$ of an $m$-dimensional smooth manifold $M$,
%embedded in $\mathbb{R}^n$ for some $m< n$, and let $\calH_m(M_t)$
%be the $m$-dimensional Hausdorff measure of the manifold $M_t$. Then
there exist constants $c_1$ and $C_1$, such that
$$\sup_{t\in[0,T]}\calH_m(M_t)\;\le\; c_1\,\calH_m(M) \;2^{C_1\,T\,
\|B^H\|_{\beta,T}^{1/\beta}},$$
 where
$\|B^H\|_{\beta,T}$
 is the
$\beta$-H\"{o}lder norm of $B^H$ (cf.\ \eqref{holdernorm}).
\end{theorem}

It is not hard to see that the H\"older norm $\|B^H\|^{1/\beta}_{\beta,T}$ grows no faster
than $O(T^{1+\epsilon})$ for any $\epsilon >0$, so that the overall rate of growth
of $\calH_m(M_t)$ given by Theorem \ref{theorem.growth.Hausdorff.1} is
$O(2^{CT^{2+\epsilon}})$. One should hope for something that was smaller, and
$H$-dependent, but current techniques fail to establish this.

Similarly, recent results of  Baudoin and  Coutin  \cite{BauCou07} seem to indicate that
correct growth rate should be $O(2^{CT^{2H}})$. These results, however, are based on
the rough path approach of \cite{Bau04,Lyon98}. While
Hairer and Ohashi in \cite{HairOha07} have proved existence
of a stationary solution of \eqref{equation.fBM.flow.candidate}, under conditions
on the vector fields $U_\gamma$ and assuming that  for $|I|<\infty$,
an approach via rough paths also seems unable to reach a better  growth rate.

Our proof of Theorem \ref{theorem.growth.Hausdorff.1} and the more general
  Theorem \ref{theorem.fBM.growth.Hausdorff.measure}
will be based on the approach  of
Hu and Nualart \cite{HuNua06}, who  obtained growth estimates on the solution
of \eqref{equation.fBM.flow.candidate}. The details follow in the remaining two
sections.

In Section \ref{section.prelim}, apart from being more formal about
setting up  notation, we shall recall some basic formulae from the
fractional calculus required for our main analysis.
Estimates on the tangent flow and the flow itself, together
with the proof of the main results, will form the bulk of Section
\ref{section.main}.

%In Section \ref{section.future}, we indicate few potentially
%promising ways of improving our results and pose few open problems,
%solutions to which are sure to improve our results.

%---------------------------------------------------------------
%---------------------------------------------------------------

\section{Preliminaries}
\label{section.prelim}

We start by listing some of the basic formulae required from the
deterministic fractional calculus, and the fractional spaces
associated with them. (See \cite{Zah98,Zah99} for a complete account
of fractional calculus.)

For $a,b\in\mathbb{R}$, $a<b$, let $L^p(a,b)$, $p\ge 1$, be the
space of Lebesgue measurable functions $f:[a,b]\to\mathbb{R}$ with
$\|f\|_{L^p(a,b)}<\iy$, where
$$\|f\|_{L^p(a,b)}= \left\{ \begin{array}{lcl}
(\int_a^b|f(x)|^p\;dx)^{\frac{1}{p}}, & \mbox{if} & 1\le p<\iy\\
ess\sup {|f(x)|:x\in [a,b]}, & \mbox{if} & p=\iy.
\end{array}\right.$$

The left sided fractional Riemann-Liouville integral of $f\in
L^1(a,b)$ of order $\al>0$ is given by
$$I^{\al}_{a+}f(x)=\frac{1}{\Gamma(\al)}\int_a^x (x-y)^{\al
-1}f(y)\;dy,$$ for almost all $x\in(a,b)$, where $\Gamma(\al)$ is
the standard Euler function. Similarly, the right sided fractional
integral is defined, for almost all $x\in(a,b)$, as
$$I^{\al}_{b-}f(x)=\frac{(-1)^{-\al}}{\Gamma(\al)}\int_x^b (y-x)^{\al
-1}f(y)\;dy,$$ where $(-1)^{-\al} =e^{-i\pi\al}$. If we consider the
fractional integral $I^{\al}_{a+}$ (or $I^{\al}_{b-}$) as an
operator with domain $L^p(a,b)$, then its range is denoted by
$I^{\al}_{a+}(L^p(a,b))$ (or $I^{\al}_{b-}(L^p(a,b))$). Clearly, for
$\al=1$, $I^{\al}_{a+}$ is the standard left integral operator, and
a simple calculation yields that $\lim_{\al\to
0}(I^{\al}_{a+}f)(x)=f(x-)=\lim_{\eps\downarrow 0}f(x-\eps)$, for
each $x\in (a,b)$. An immediate consequence of the definition of the
fractional integral is that
\begin{equation}
\label{equation.fBM.composition.left.fractional.integral}
I^{\al}_{a+}(I^{\beta}_{a+}f)=I^{\al+\beta}_{a+}f,
\end{equation}
for all $\al,\beta >0$. With some obvious variations where
needed, all  hold also for right sided
fractional integrals; viz.
\begin{eqnarray*}
(I^1_{b-}f)(x) & = & (-1)\int_x^b f(y)\;dy, \\
\lim_{\al\to 0} (I^{\al}_{b-}f)(x) & = & f(x+)
\ = \ \lim_{\eps\downarrow 0}f(x+\eps),\\
I^{\al}_{b-}(I^{\beta}_{b-}f)& = & I^{\al+\beta}_{b-}f,\; \forall
\al,\beta >0
\end{eqnarray*}
(See \cite{Zah98} for these and more on fractional calculus.)

Having defined a fractional integral, we now define a fractional
derivative as the  inverse of the fractional integral operator,
whenever it is well defined. In other words, to each element $f$ in
$I^{\al}_{a+}(L^p(a,b))$ there corresponds a $\phi\in L^p(a,b)$,
such that $I^{\al}_{a+}\phi=f$. This $\phi$ is unique in $L^p(a,b)$ and
agrees almost everywhere with the fractional derivative,
known as the  left sided Riemann-Liouville or Weyl derivative,
of $\al^{th}$-order and  defined as
\begin{eqnarray}
\label{equation.fBM.left.frac.derivative} D^{\al}_{a+}f(x) & = &
\Big(\frac{1}{\Gamma(1-\al)}\frac{d}{dx}\int_a^x
\frac{f(y)}{(x-y)^{\al}}\;dy\Big)\;1_{(a,b)}(x)\nonumber\\
& = & \frac{1}{\Gamma(1-\al)}\Big(\frac{f(x)}{(x-a)^{\al}}
+\al\int_a^x \frac{f(x)-f(y)}{(x-y)^{1+\al}}\;dy\Big)\;1_{(a,b)}(x).
\end{eqnarray}
Equivalently, we can write $D^{\al}_{a+}f=D(I^{1-\al}_{a+}f),$ where
$D$ is the standard derivative operator. Similarly, we can define
the right sided Weyl derivative as
$D^{\al}_{b-}f=D(I^{1-\al}_{b-}f)$, for which
\begin{eqnarray}
\label{equation.fBM.right.frac.derivative} D^{\al}_{b-}f(x) & = &
\Big(\frac{(-1)^{\al-1}}{\Gamma(1-\al)}\frac{d}{dx}\int_x^b
\frac{f(y)}{(y-x)^{\al}}\;dy\Big)\;1_{(a,b)}(x)\nonumber\\
& = & \frac{(-1)^{\al}}{\Gamma(1-\al)}\Big(\frac{f(x)}{(b-x)^{\al}}
+\al\int_x^b \frac{f(x)-f(y)}{(y-x)^{1+\al}}\;dy\Big)\;1_{(a,b)}(x).
\end{eqnarray}

As in the case of the integral operators, there is an analogue of
the composition formula, given, for all $\al, \beta >0$, by
\begin{equation}
\label{equation.fBM.composition.left.fractional.derivative}
D^{\al}_{a+}(D^{\beta}_{a+}f)=D^{\al+\beta}_{a+}f.
\end{equation}
A similar formula also holds for the right sided derivatives, and is
given by,
\begin{equation}
\label{equation.fBM.composition.right.fractional.derivative}
D^{\al}_{b-}(D^{\beta}_{b-}f)=D^{\al+\beta}_{b-}f,
\end{equation}
as long as all the fractional derivatives are well defined.

We note that the linear spaces $I^{\al}_{a+}(L^p(a,b))$, for various
choices of $\al$ and $p$, are Banach spaces equipped with the norms
$$\|f\|_{I^{\al}_{a+}(L^p(a,b))}= \|f\|_{L^p(a,b)}+
\|D^{\al}_{a+}f\|_{L^p(a,b)},$$ and a similar norm is defined on the
space $I^{\al}_{b-}(L^p(a,b))$.

%Using the methods of fractional calculus, one can extend the
%standard integration by parts formula to the more general case of
%$L^p$ functions. Hence the generalized {\em integration by parts
%formula} can be written as \be\label{equation.fBM.int.by.parts.1}
%\int_a^b f(x)I^{\al}_{a+}g(x)dx=(-1)^{\al}\int_a^b
%g(x)I^{\al}_{b-}f(x)dx, \ee where $f\in L^p(a,b),\; g\in L^q,\; p\ge
%1,\; q\ge 1,\; 1/p+1/q\le 1+\al,\; p,q>1,$ and $1/p+1/q= 1+\al$. A
%similar formula, called the {\em second integration by parts
%formula}, holds true for derivative operators, and is given as \be
%\label{equation.fBM.int.by.parts.2}(-1)^{\al}\int_a^b
%D^{\al}_{a+}f(x) g(x)dx=\int_a^b f(x)D^{\al}_{b-}g(x)dx, \ee where
%$f\in I^{\al}_{a+}(L^p(a,b)),\;g\in I^{\al}_{b-}(L^q(a,b)),\;
%p\ge1,\;q\ge1,\;1/p+1/q\le 1+\al.$

Let $f(a+)=\lim_{\eps\downarrow 0}f(a+\eps),$ and $g(b-)=
\lim_{\eps\downarrow 0}f(b-\eps)$, whenever the limit exists and is
finite, and define
\begin{eqnarray*}
f_{a+}(x) & = & (f(x)-f(a+))1_{(a,b)}(x),\\
g_{b-}(x) & = & (g(x)-g(b-))1_{(a,b)}(x).
\end{eqnarray*}
Using the methods of fractional calculus (see \cite{Zah98}), an
extension of the Stieltjes integral, called the {\em generalized
Stieltjes integral}, of $f$ with respect to $g$ can be defined as
\begin{equation}
\label{equation.fBM.define.int.2}
\int_a^bf(x)dg(x)=(-1)^{\al}\int_a^b D^{\al}_{a+}f(x)
D^{1-\al}_{b-}g_{b-}(x)\,dx,
\end{equation}
where $f_{a+}\in
I^{\al}_{a+}(L^p(a,b))$ and $g_{b-}\in I^{1-\al}_{b-}(L^q(a,b))$ for
some $p,q\ge 1,\;1/p+1/q\le 1,\; 0\le\al\le 1,$ and $\al p<1$.

Next we define $C^{\lambda}(a,b; \mathbb{R}^d)$, the space of
$\lambda$-H\"{o}lder continuous functions, with $\lambda\in (0,1]$,
as the space of $\mathbb{R}^d$ valued functions for some fixed $d\in
\mathbb{N}$, the set of natural numbers, equipped with the norm
given by
\begin{eqnarray}
\label{holdernorm}
 \|f\|_{\lambda,a,b} := \sup_{a\le c\le
d\le b}\frac{\|f(d)-f(c)\|_2}{|d-c|^{\lambda}},
\end{eqnarray}
where $\|\cdot\|_2$ is the usual Euclidean norm in the appropriate dimension.
(When $a=0$, we shall write $ \|f\|_{\lambda,b}$ for $\|f\|_{\lambda,0,b}$.)

In \cite{Zah98}, Z\"{a}hle proved that the conditions of the
definition \eqref{equation.fBM.define.int.2} are met if $f\in
C^{\lambda}(0,T; \mathbb{R})$ and $g\in C^{\mu}(0,T; \mathbb{R})$
for $\lambda +\mu
>1$, in which case the integral defined in \eqref{equation.fBM.define.int.2}
coincides with the Riemann-Stieltjes integral. Now we state the
following well known result concerning the H\"{o}lder coefficient
and exponent of fractional Brownian motion with Hurst parameter $H$.

\begin{lemma}
\label{lemma.fBM.Holder} For $\{B^H(t):\; t\in [0,T]\}$, a
fractional Brownian motion with Hurst parameter $H\in (0,1)$, there
exists,  for each $0<\eps<H$ and $T>0$,
 a positive random variable $\eta_{\eps,T}$, such that $E(|\eta_{\eps,T}|^p)<\iy$ for all
$p\in[1,\iy)$ and, for all $s,t\in [0,T]$,
$$|B^H(t)-B^H(s)|\;\le\; \eta_{\eps,T}|t-s|^{H-\eps}\;\;\mbox{a.s.,}$$
where $\eta_{\eps,T} = C_{H,\eps}T^{H-\eps}\xi_T$, with the
$L^q(\Omega)$ norm of $\xi_T$ bounded by $c_{\eps,q}T^{\eps}$ for
$q\ge \frac{2}{\eps}$.
\end{lemma}

(For a proof of this, which involves a simple application of
a Garsia-Rodemich-Rumsey type inequality, we refer the reader to
\cite{NuaRas02}.)

Consequently, using the above theory of deterministic fractional integration,
 integrals with respect to the fractional Brownian motion
can also be defined, for appropriate integrands. This was done in
\cite{Zah98}, where a corresponding stochastic calculus is also developed
with an appropriate  change of variables formula.

For the following definitions, we shall assume $\al<\frac{1}{2}$.

Define $W^{1-\al, \infty}_T(0,T; \mathbb{R})$ to be the space of
measurable functions $g:[0,T]\to\mathbb{R}$, endowed with and finite under the
 norm
\be\label{equation.fBM.W.1.alpha.infinity.T}
\|g\|_{1-\al,\iy,T}:=\sup_{0<s<t<T}
\big(\frac{|g(t)-g(s)|}{(t-s)^{1-\al}}
+\int_s^t\frac{|g(y)-g(s)|}{|t-s|^{2-\al}} dy \big). \ee
Clearly, \be\label{equation.fBM.embedding.C.in.W.2}
C^{1-\al+\eps}(0,T; \mathbb{R})\,\subset\, W^{1-\al,\iy}_T(0,T;
\mathbb{R})\,\subset\, C^{1-\al}(0,T; \mathbb{R}), \ee for all
$\eps>0$. Moreover, if $g\in W^{1-\al, \infty}_T(0,T; \mathbb{R})$,
then $g_{|_{(0,t)}}\in I^{1-\al}_{t-}(L^{\iy}(0,t))$ for all $t\in
(0,T)$.

Recalling now the vector fields introduced in
\eqref{equation.fBM.flow.candidate}, the time has come to demand a set of
regularity assumptions. Assume that there exist constants $M_{\gamma}$,
$M_{\gamma}^{(1)}$ and $M_{\gamma}^{(2)}$ for all $\gamma\in\mathbb{N}$
such that:

\bi
\item[$(A1)$]
$|U^i_{\gamma}(x)| \le M_{\gamma}$, $\forall x\in \mathbb{R}^n$ and
$\gamma\in\mathbb{N}$, where $U^i_{\gamma}$ denotes the $i$-th
component of $U_\gamma$.

\item[$(A2)$]
$|U^i_{\gamma}(x)-U^i_{\gamma}(y)|
\le M_{\gamma}^{(1)}\|x-y\|_2$, $\forall x,y\in \mathbb{R}^n$ and
$\gamma\in\mathbb{N}$, where $\|\cdot\|_2$ denotes the standard Euclidean
norm in the appropriate dimension.

\item[$(A3)$]
$|W_{\gamma,\,j}^{i}(x)-W_{\gamma,\,j}^{i}(y)|\le
M^{(2)}_{\gamma}\|x-y\|_2$, $\forall x,y\in \mathbb{R}^n$ and
$\gamma\in\mathbb{N}$, where $W_{\gamma}(x)$ denotes the spatial derivative of
$U_{\gamma}(x)$, and $W_{\gamma,\,j}^{i}(\cdot)$ denotes the
$(i,j)$-th element of the matrix $W_{\gamma}(\cdot)$.

\item[$(A4)$]
$M^{(1)}=\sum_{\al\in \mathbb{N}}M^{(1)}_{\al} <\infty$,
$M^{(2)}=\sum_{\al\in \mathbb{N}}M^{(2)}_{\al} <\infty$, and
$M^{(3)}=\sum_{\al\in \mathbb{N}}M^{(3)}_{\al} <\infty$.

\ei

Under conditions $(A1)-(A4)$, existence and
uniqueness of the solution of \eqref{equation.fBM.flow.candidate},
in the space $C^{1-\al}(0,T;\mathbb{R}^n)$ is proven in
\cite{NuaRas02} for $|I|<\iy$.

In fact, the existence and uniqueness of the solution can be proven
under far weaker conditions, but without necessarily giving a solution which
provids a diffeomorphism in $\mathbb{R}^n$ (cf.\  \cite{DecNua06} for
details). Properties  of the solution of the flow equation
are also obtained in \cite{NuaRas02}, and improved on in
\cite{HuNua06}.

%------------------------------------------------------------------------------------------
%------------------------------------------------------------------------------------------

\section{The main result}
\label{section.main}

We shall now adopt and adapt the approach developed in \cite{HuNua06} to derive some
estimates on some of the basic geometric characteristics of the flow
 \eqref{equation.fBM.flow.candidate}.

%We start with recalling the definition of the flow in
%\eqref{equation.fBM.flow.candidate} and class of vector fields
%$\{U_{\gamma}\}_{\gamma\in\mathbb{N}}$ defined in
%\eqref{equation.fBM.vector.fields}, together with the corresponding
%Assumptions $(A1)$-$(A4)$.

With $M$, as usual, a $C^2$, $m$-dimensional  manifold embedded in
 $\mathbb{R}^n$,  we write $T_xM$ for its tangent
space at $x$. Let $v\in T_xM$. Then its push-forward
under the flow $\Phi_t$ is denoted by
$$v_t=D\Phi_t(x)v,$$
where $D\Phi_t(x)=(\frac{\partial\Phi_t^i(x)}{\partial x^j})_{ij}$
denotes the matrix of spatial derivatives of the flow $\Phi_t(x)$,
and $v_t\in T_{x_t}M_t$. From  now on we shall write $x_t$ for
$\Phi_t(x)$.

We now
prove the following technical result, which will form the basis for
much of the subsequent analysis.
%Before we shall present the main
%result of this part of the thesis, we shall state and prove a
%technical result which will lead us to the main result of Part
%\ref{part.fractional}.

\begin{theorem} \label{theorem.fBM.growth.vector} Under assumptions
$(A1)-(A4)$, and for $\al= 1-H+\delta$, $\beta=H-\eps$,
such that $(1-H)<\al<1/2$ and $\delta>\eps$, there exist a constant
$c$ and a random variable $C_T$, such that
\begin{eqnarray*}
\sup_{r\in[0,T]}\|v_r\|_2 & \le & \sup_{r\in[0,T]}\|v_r\|_1\\
& \le & c\;2^{C_T\;T},
\end{eqnarray*}
where $\|v_r\|_2$ and $\|v_r\|_1$ denote the $l_2$ and $l_1$ norms,
respectively, of the vector $v_r$ as an element in $\mathbb{R}^n$.
The random variable $C_T$ depends on $\al$, $\beta$, $n$, $I$, and
$\{\|B^H_{\gamma}\|_{\beta,T}, M_{\gamma}, M^{(1)}_{\gamma},
M^{(2)}_{\gamma}\}_{\gamma\in I}$. Furthermore,
$$E[C_T]^{\beta} \le C\cdot E [\|B^H\|_{\beta,T}],$$ where the constant $C$
depends only on $\al$, $\beta$, $n$, $|I|$ and $\{M_{\gamma},
M^{(1)}_{\gamma}, M^{(2)}_{\gamma}\}_{\gamma\in I}$.
\end{theorem}

\begin{remark}
\label{remark.fBM.growth.vector.1} For a better understanding of the
results of Theorem \ref{theorem.fBM.growth.vector}, we note that for
the case $|I|=1$, this result simplifies to
\begin{eqnarray*}
\sup_{r\in[0,T]}\|v_r\|_2 & \le &
c\;2^{C\;T\;\|B^H\|^{1/\beta}_{\beta,T}},
\end{eqnarray*}
for some constants $c$ and $C$, dependent only on the various
uniform bounds and the Lipschitz coefficients corresponding to the
vector field.
\end{remark}

\begin{remark}
\label{remark.fBM.growth.vector.2} The results listed in this
section hold true for any $I\subset \mathbb{N}$ as long as the
cardinality of the set satisfies $|I|<\iy$. However, extensions of
these results to the case $I=\mathbb{N}$, though possible, require
unnatural conditions on the summability of the constants appearing
in Assumptions $(A1)-(A3)$. For instance, extending Lemma
\ref{lemma.fBM.flow.est} to the case $I=\mathbb{N}$ would require
$$\frac{\sum_{\gamma\in \mathbb{N}}M^{(1)}_{\gamma}\|B^H_{\gamma}\|_{\beta,T}}{
\sum_{\gamma\in
\mathbb{N}}M^{(2)}_{\gamma}\|B^H_{\gamma}\|_{\beta,T}} <\iy.$$
This, in turn would be implied by
$\sum_{\gamma\in \mathbb{N}} {M^{(1)}_{\gamma}}/{M^{(2)}_{\gamma}} <\iy$,
which does not seem to have a clear meaning in terms of
the vector fields $U_{\gamma}$.
\end{remark}

The idea of the proof of Theorem \ref{theorem.fBM.growth.vector} is
to break up the interval $[0,T]$ into smaller units of size
$\Delta$, on which reasonable estimates of $\|v_r\|_2$ are possible,
and then to glue the intervals together to obtain the required
result. However, in the process, we shall need to derive an
estimate on the flow,  presented in the following lemma, the proof of which
relies on some results of \cite{HuNua06}.

\begin{lemma}
\label{lemma.fBM.flow.est} Let $M$, $M^{(1)}$ be constants as
defined in Assumptions $(A1)-(A4)$, and $0\le s\le t\le T$ be such
that
$$(t-s)^{-\beta}>\frac{n\,\al(2\al+\beta-1)}{2(1-\al)(1-2\al)(\al+\beta-1)\Gamma(\al)\Gamma(1-\al)}
\sum_{\gamma\in I} M^{(1)}_{\gamma}\|B^H_{\gamma}\|_{\beta,T},$$
where $\al= 1-H+\delta$, $\beta=H-\eps$, such that $(1-H)<\al<1/2$
and $\delta>\eps$. Then for $x_t$ defined in
\eqref{equation.fBM.flow.candidate} there exists a positive random variable
$K^*_{s,t}$ such that
\begin{equation}
\label{equation.fBM.flow.est} \int_s^t
\frac{\|x_t-x_r\|_2}{(t-r)^{1+\al}}\;dr\; \le\;
K^*_{s,t}(t-s)^{\beta-\al}.
\end{equation}
Furthermore, $K^*_{s,t}$ can be bounded above by another random
variable, independent of $s$ and $t$, with finite moments of order
greater than $1$, as long as $(t-s)$ is chosen sufficiently small.
\end{lemma}

\begin{remark}
\label{remark.fBM.alpha.beta} Note that under the aforementioned
conditions concerning $\al$ and $\beta$, we have $\al+\beta>1$, and
$\beta>\al$.
\end{remark}

\textbf{Proof:} Writing $U^i_{\gamma}(\cdot)$ for the $i$-th
component of the vector $U_{\gamma}(\cdot)$ and choosing
$\{e_i\}_{i=1}^n$ as the canonical basis of $\mathbb{R}^n$, we have
$$\langle (x_t-x_s),e_i\rangle = \sum_{\gamma\in I}\int_s^t U^i_{\gamma}(x_r)\;dB^H_{\gamma}(r),$$
which is true by linearity of the operation, and where
$\langle\cdot,\cdot\rangle$ denotes the standard Euclidean inner
product. Hence for $\al\in (1-H,\frac{1}{2})$, using \eqref{equation.fBM.define.int.2},
we obtain
\begin{eqnarray*}
|\langle (x_t-x_s),e_i\rangle| & = & \Big|\sum_{\gamma\in I}\int_s^t
U^i_{\gamma}(x_r)\;dB^H_{\gamma}(r)\Big|\\
& = & \Big|\sum_{\gamma\in I}\int_s^t D^{\al}_{s+}U^i_{\gamma}(x_r)
D^{1-\al}_{t-}B_{\gamma,t-}^H(r)\;dr\Big|\\
& \le & \sum_{\gamma\in I}\int_s^t
|D^{\al}_{s+}U^i_{\gamma}(x_r)|\cdot
|D^{1-\al}_{t-}B_{\gamma,t-}^H(r)|\;dr
\end{eqnarray*}
To obtain a bound on the second term in the integrand, choose $\beta
<H$, such that $\al+\beta>1$, so that using
\eqref{equation.fBM.right.frac.derivative}, we have
\begin{eqnarray}\label{equation.fBM.est.1}
|D^{1-\al}_{t-}B^H_{\gamma,t-}(r)| & = &
\Big|\frac{(-1)^{1-\al}}{\Gamma(\al)} \Big(
\frac{B^H_{\gamma}(t)-B^H_{\gamma}(r)}{(t-r)^{1-\al}}+
\al\int_r^t\frac{B^H_{\gamma}(u)-B^H_{\gamma}(r)}{(u-r)^{2-\al}}du\Big)\Big|\nonumber\\
& \le &
\frac{1}{\Gamma(\al)}\Big(\frac{|B^H_{\gamma}(t)-B^H_{\gamma}(r)|}{|t-r|^{1-\al}}
+\al\int_r^t\frac{|B^H_{\gamma}(u)-B^H_{\gamma}(r)|}{(u-r)^{2-\al}}du\Big)\nonumber\\
& = &
\frac{1}{\Gamma(\al)}\Big(\frac{|B^H_{\gamma}(t)-B^H_{\gamma}(r)|
(t-r)^{\beta}}{(t-r)^{\beta}(t-r)^{1-\al}}\nonumber\\
&   & \mbox{}+
\al\int_r^t\frac{|B^H_{\gamma}(u)-B^H_{\gamma}(r)|}{(u-r)^{\beta}}
(u-r)^{\al+\beta-2}du\Big)\nonumber\\
& \le & \frac{1}{\Gamma(\al)}
\Big(\|B^H_{\gamma}\|_{\beta,T}(t-r)^{\al+\beta-1}+\al\|B^H_{\gamma}\|_{\beta,T}
\frac{(t-r)^{\al+\beta-1}}{\al+\beta-1}\Big)\nonumber\\
& = & k_1(\al,\beta)\|B^H_{\gamma}\|_{\beta,T} (t-r)^{\al+\beta-1},
\end{eqnarray}
where $k_1(\al,\beta) =
\frac{(2\al+\beta-1)}{(\al+\beta-1)\Gamma(\al)}$.

To bound the first term we use
\eqref{equation.fBM.left.frac.derivative} and assumptions
$(A1)-(A2)$ to see that
\begin{eqnarray}\label{equation.fBM.est.2}
|D^{\al}_{s+}U^i_{\gamma}(x_r)| & = & \frac{1}{\Gamma(1-\al)}
\Big|\frac{U^i_{\gamma}(x_r)}{(r-s)^{\al}}+\al\int_s^r
\frac{(U^i_{\gamma}(x_r)-U^i_{\gamma}(x_{\theta}))}{(r-\theta)^{1+\al}}\;d\theta\Big|\nonumber \\
& \le & \frac{1}{\Gamma(1-\al)}
\Big(\frac{|U^i_{\gamma}(x_r)|}{(r-s)^{\al}}+\al\int_s^r
\frac{|U^i_{\gamma}(x_r)-U^i_{\gamma}(x_{\theta})|}{(r-\theta)^{1+\al}}\;d\theta\Big)\nonumber \\
& \le & c_{\al} \Big(\frac{M_{\gamma}}{(r-s)^{\al}}+\al\int_s^r
\frac{M^{(1)}_{\gamma}\|x_r-x_{\theta}\|_2}{(r-\theta)^{1+\al}}\;d\theta\Big)\nonumber \\
& \le & c_{\al} \Big(M_{\gamma}(r-s)^{-\al}+
M^{(1)}_{\gamma,\al}\|x\|_{s,r,1-\al}(r-s)^{1-2\al}\Big),
\end{eqnarray}
where $c_{\al}=\Gamma(1-\al)^{-1}$, $M^{(1)}_{\gamma,\al}=
\frac{\al M^{(1)}_{\gamma}}{(1-2\al)}$ and $\|x\|_{s,r,1-\al}$ is the
H\"{o}lder norm as defined in \eqref{holdernorm}.

Therefore, combining the above two estimates , we find
\begin{eqnarray*}
|\langle (x_t-x_s),e_i\rangle| & \le &
c_{\al}k_1(\al,\beta)\sum_{\gamma\in I}\|B^H_{\gamma}\|_{\beta,T}
\int_s^t\Big(M_{\gamma}(r-s)^{-\al}(t-r)^{\al+\beta-1}\\
&    & \mbox{}+ M^{(1)}_{\gamma,\al} \|x\|_{s,r,1-\al}
(r-s)^{1-2\al}(t-r)^{\al+\beta-1}\Big)\;dr\\
& \le & c_{\al}k_1(\al,\beta) \sum_{\gamma\in
I}\|B^H_{\gamma}\|_{\beta,T} (t-s)^{\al+\beta-1}
\int_s^t\Big(M_{\gamma}(r-s)^{-\al} \\
&    &\mbox{}+ M^{(1)}_{\gamma,\al} \|x\|_{s,r,1-\al}
(r-s)^{1-2\al}\Big)\;dr\\
& \le & c_{\al}k_1(\al,\beta)\sum_{\gamma\in I} \|B^H_{\gamma}\|_{\beta,T}\Big(M_{\gamma}(t-s)^{\beta}(1-\al)^{-1}\\
&   & \mbox{}+ M^{(1)}_{\gamma,\al}
\|x\|_{s,t,1-\al}(t-s)^{1-\al+\beta}(2-2\al)^{-1}\Big).
\end{eqnarray*}

Let
\begin{equation}
\label{equation.fBM.constt.M(alpha)} M_{\al}=
(1-\al)^{-1}\sum_{\gamma\in I}M_{\gamma}\|B^H_{\gamma}\|_{\beta,T},
\end{equation}
and
\begin{equation}
\label{equation.fBM.constt.M(tilde.1.alpha)} \tilde{M}^{(1)}_{\al}=
(2-2\al)^{-1} \sum_{\gamma\in
I}M^{(1)}_{\gamma,\al}\|B^H_{\gamma}\|_{\beta,T}.
\end{equation}
Then
\begin{eqnarray*}
\|(x_t-x_s)\|_1 & = & \sum_{i=1}^n
|\langle(x_t-x_s),e_i\rangle|\nonumber\\
& \le & c_{\al}n k_1(\al,\beta) \Big(M_{\al}(t-s)^{\beta} +
\tilde{M}^{(1)}_{\al} \|x\|_{s,t,1-\al}(t-s)^{1-\al+\beta}\Big).
\end{eqnarray*}

Equivalently,
\begin{eqnarray}
\label{equation.fBM.flow.est.1}
\frac{\|(x_t-x_s)\|_1}{(t-s)^{1-\al}} & \le &
c_{\al}nk_1(\al,\beta)\Big(M_{\al}(t-s)^{\al+\beta-1}
\nonumber\\
&    & \mbox{}+ \tilde{M}^{(1)}_{\al}
\|x\|_{s,t,1-\al}(t-s)^{\beta}\Big).
\end{eqnarray}
(Recall that $\al+\beta>1$.)

Now using the above estimate, and the fact that $\|\cdot\|_2$ is
bounded above by $\|\cdot\|_1$, we have
\begin{eqnarray}\label{equation.fBM.flow.est.2}
\|x\|_{s,t,1-\al} & = & \sup_{s\le u\le v\le
t}\frac{\|(x_v-x_u)\|_2}{(v-u)^{1-\al}}\nonumber\\
& \le & \sup_{s\le u\le v\le
t}\frac{\|(x_v-x_u)\|_1}{(v-u)^{1-\al}}\nonumber\\
& \le & \sup_{s\le u\le v\le t}
c_{\al}nk_1(\al,\beta)\Big(M_{\al}(v-u)^{\al+\beta-1}\nonumber\\
&    & \mbox{}+
\tilde{M}^{(1)}_{\al} \|x\|_{u,v,1-\al}(v-u)^{\beta}\Big)\nonumber\\
& \le &
c_{\al}nk_1(\al,\beta)\Big(M_{\al}(t-s)^{\al+\beta-1}\nonumber\\
&    & \mbox{}+ \tilde{M}^{(1)}_{\al}
\|x\|_{s,t,1-\al}(t-s)^{\beta}\Big).
\end{eqnarray}

Now choosing $s,t$ such that
\begin{eqnarray}
\label{equation.fBM.interval.size.1} (t-s)^{-\beta} & > &
c_{\al}nk_1(\al,\beta)\tilde{M}^{(1)}_{\al},
\end{eqnarray}
\eqref{equation.fBM.flow.est.2} can be rewritten as
\begin{eqnarray}
\label{equation.fBM.flow.pre.est} \|x\|_{s,t,1-\al} & \le &
\frac{c_{\al}nk_1(\al,\beta)M_{\al}(t-s)^{\al+\beta-1}}{1-
c_{\al}nk_1(\al,\beta)\tilde{M}^{(1)}_{\al}(t-s)^{\beta}}\nonumber\\
& = & K_{s,t} (t-s)^{\al+\beta-1},
\end{eqnarray}
where $K_{s,t}=\frac{c_{\al}nk_1(\al,\beta)M_{\al}}{1-
c_{\al}nk_1(\al,\beta)\tilde{M}^{(1)}_{\al}(t-s)^{\beta}}$.

Therefore,
\begin{eqnarray*}
\int_s^t \frac{\|x_t-x_r\|_2}{(t-r)^{1+\al}}\;dr & = & \int_s^t
\frac{\|x_t-x_r\|_2}{(t-r)^{1-\al}}(t-r)^{-2\al}\;dr\\
& \le & \|x\|_{s,t,1-\al}\int_s^t (t-r)^{-2\al}\;dr\\
& \le & K_{s,t}\frac{(t-s)^{\beta-\al}}{(1-2\al)}\\
& = & K^*_{s,t}(t-s)^{\beta-\al},
%& \le & \frac{M_{\al}}{(c-1)\tilde{M}^{(1)}_{\al}}
%\cdot\frac{1}{(t-s)^{1-\al}}\cdot\frac{(t-s)^{1-2\al}}{(1-2\al)}\\
%& = & \frac{M_{\al}}{(c-1)(1-2\al)\tilde{M}^{(1)}_{\al}}\cdot
%\frac{1}{(t-s)^{\al}}\\
%& = & \frac{2\sum_{\gamma\in
%I}M_{\gamma}\|B^H_{\gamma}\|_{\beta,T}}{\al\,(c-1)\sum_{\gamma\in I}
%M^{(1)}_{\gamma}\|B^H_{\gamma}\|_{\beta,T}}\cdot
%\frac{1}{(t-s)^{\al}}\\
%& \le & \frac{2}{\al\,(c-1)}\sum_{\gamma\in
%I}\frac{M_{\gamma}}{M^{(1)}_{\gamma}}\cdot
%\frac{1}{(t-s)^{\al}}\\
%& = & \frac{K^*}{(t-s)^{\al}},
\end{eqnarray*}
where $K^*_{s,t}=\frac{K_{s,t}}{(1-2\al)},$ thus establishing
\eqref{equation.fBM.flow.est}. The final claim, that $K^*_{s,t}$ can
be bounded by a random variable independent of $s$ and $t$, will be
proven later. \\

%Note that a similar bound works for all intervals of size smaller
%than $(t-s)$. In particular, the random coefficient $K^*_{s,t}$
%decreases as the size of the interval decreases.

\textbf{Proof of Theorem {\ref{theorem.fBM.growth.vector}}:}
Taking the space derivative of \eqref{equation.fBM.flow.candidate},
the existence of which is ensured by Theorem 3.2 in \cite{HuNua06},
we have
$$D\Phi_t(x)=I +\sum_{\gamma\in I}\int_0^t W_{\gamma}(\Phi_s(x))D\Phi_s(x)\;dB^H_{\gamma}(s),$$
where the matrix $W_{\gamma}(\cdot)=(W^i_{\gamma,j}(\cdot))_{i,j}$
denotes the spatial derivative of the vector field $U$.

Now using the definition of the pushforward of a vector, we can
write the evolution equation of the tangent vector as follows
$$v_t = v + \sum_{\gamma\in I}\int_0^t W_{\gamma}(x_s)v_s dB^H_{\gamma}(s).$$

Recall that $\|v_t\|_1 =\sum_{i=1}^n |\langle v_t,e_i \rangle|$,
where $\langle\cdot,\cdot\rangle$ is the standard Euclidean inner
product, and $\{e_i\}_{i=1}^n$ denotes the canonical basis of
$\mathbb{R}^n$. Since,
$$\langle v_t,e_i\rangle = x+ \sum_{\gamma\in I}\int_0^t
\langle W_{\gamma}(x_r)v_r,e_i\rangle dB^H_{\gamma}(r),$$ we have
\begin{eqnarray*}
|\langle v_t,e_i\rangle-\langle v_s,e_i\rangle| & = & \Big|
\sum_{\gamma\in I}\int_s^t\langle W_{\gamma}(x_r)v_r,e_i\rangle dB^H_{\gamma}(r)\Big|\\
& = & \Big|\sum_{\gamma\in I}\int_s^t D^{\al}_{s+}\langle
W_{\gamma}(x_r)v_r,e_i\rangle
D^{1-\al}_{t-}B^H_{\gamma,t-}(r) dr\Big|\\
& \le & \sum_{\gamma\in I}\int_s^t |D^{\al}_{s+}\langle
W_{\gamma}(x_r)v_r,e_i\rangle| \cdot
|D^{1-\al}_{t-}B^H_{\gamma,t-}(r)| dr.
\end{eqnarray*}
The above inequality holds for any choice of $s$ and $t$, but we are
interested in pairs for which $(t-s)$ is sufficiently small. To this
end, note first that from \eqref{equation.fBM.est.1} we can bound
the second integrand by
$$ |D^{1-\al}_{t-}B^H_{\gamma,t-}(r)| \le k_1(\al,\beta)\|B^H_{\gamma}\|_{\beta,T}
(t-r)^{\al+\beta-1}.$$

Now using \eqref{equation.fBM.left.frac.derivative} and Assumptions
(A2)--(A4), the first integrand can be bounded by
\begin{eqnarray*}
|D^{\al}_{s+}\langle W_{\gamma}(x_r)v_r,e_i\rangle| & \le &
\frac{1}{\Gamma(1-\al)}\Big[\frac{|\langle
W_{\gamma}(x_r)v_r,e_i\rangle|}{(r-s)^{\al}}
+\al\int_s^r\frac{|\langle W_{\gamma}(x_r)v_r,e_i\rangle-\langle
W_{\gamma}(x_{\theta})v_{\theta},e_i\rangle|}{(r-\theta)^{1+\al}}d\theta
\Big]\\
& \le & \frac{1}{\Gamma(1-\al)}\Big[\frac{\sum_{j=1}^n|
W^i_{\gamma,j}(x_r)\langle v_r,e_j\rangle|}{(r-s)^{\al}}
+\al\int_s^r\frac{|\langle W_{\gamma}(x_r)v_r,e_i\rangle-\langle
W_{\gamma}(x_{\theta})v_{\theta},e_i\rangle|}{(r-\theta)^{1+\al}}d\theta
\Big]\\
& \le &
\frac{1}{\Gamma(1-\al)}\Big[M^{(1)}_{\gamma}\sum_{j=1}^n\frac{|\langle
v_r,e_j\rangle|}{(r-s)^{\al}} +\al\int_s^r\frac{|\langle
W_{\gamma}(x_r) v_r,e_i\rangle- \langle W_{\gamma}(x_{\theta})
v_r,e_i\rangle|}{(r-\theta)^{1+\al}}d\theta\\
&     & \mbox{}+ \al\int_s^r\frac{|\langle W_{\gamma}(x_{\theta})
v_r,e_i\rangle- \langle
W_{\gamma}(x_{\theta})v_{\theta},e_i\rangle|}{(r-\theta)^{1+\al}}d\theta\Big]\\
& = & \frac{1}{\Gamma(1-\al)}
\Big[M^{(1)}_{\gamma}\sum_{j=1}^n\frac{|\langle
v_r,e_j\rangle|}{(r-s)^{\al}}\\
&   & \mbox{} +\al\int_s^r\frac{|\sum_{j=1}^n
(W^i_{\gamma,j}(x_r)\langle v_r,e_j\rangle-
W^i_{\gamma,j}(x_{\theta})\langle
v_r,e_j\rangle)|}{(r-\theta)^{1+\al}}d\theta\\
&     & \mbox{}+ \al\int_s^r\frac{|\sum_{j=1}^n
(W^i_{\gamma,j}(x_{\theta})\langle v_r,e_j\rangle-
W^i_{\gamma,j}(x_{\theta})\langle v_{\theta},e_i\rangle)|}{(r-\theta)^{1+\al}}d\theta\Big]\\
& \le & \frac{1}{\Gamma(1-\al)}
\Big[M^{(1)}_{\gamma}\sum_{j=1}^n\frac{|\langle
v_r,e_j\rangle|}{(r-s)^{\al}}\\
&   & \mbox{} +\al\int_s^r\frac{\sum_{j=1}^n |W^i_{\gamma,j}(x_r)-
W^i_{\gamma,j}(x_{\theta})|\cdot |\langle
v_r,e_j\rangle|}{(r-\theta)^{1+\al}}d\theta\\
&     & \mbox{}+ \al\int_s^r\frac{\sum_{j=1}^n
|W^i_{\gamma,j}(x_{\theta})|\cdot |\langle v_r,e_j\rangle-
\langle v_{\theta},e_i\rangle|}{(r-\theta)^{1+\al}}d\theta\Big]\\
& \le & \frac{1}{\Gamma(1-\al)}
\Big[M^{(1)}_{\gamma}\sum_{j=1}^n\frac{|\langle
v_r,e_j\rangle|}{(r-s)^{\al}} +\al M^{(2)}_{\gamma}\sum_{j=1}^n
|\langle v_r,e_j\rangle| \int_s^r \frac{\|x_r-
x_{\theta}\|_2}{(r-\theta)^{1+\al}}d\theta \\
&     & \mbox{}+ \al M^{(1)}_{\gamma}\sum_{j=1}^n \int_s^r
\frac{|\langle v_r,e_j\rangle- \langle
v_{\theta},e_i\rangle|}{(r-\theta)^{1+\al}}d\theta\Big].
\end{eqnarray*}

Now using the result proven in Lemma \ref{lemma.fBM.flow.est}, for
$r$ such that $s<r<t$, with $(t-s)$ satisfying
\eqref{equation.fBM.interval.size.1}, we have
\begin{equation}
\int_s^r \frac{\|x_r-x_{\theta}\|_2}{(r-\theta)^{1+\al}} \,d\theta
\le K^*_{s,r}(r-s)^{\beta-\al}\nonumber.
\end{equation}

Hence,
\begin{eqnarray*}
|D^{\al}_{s+}\langle W_{\gamma}(x_r)v_r,e_i\rangle| & \le &
\sum_{j=1}^n\Big[\frac{|\langle
v_r,e_j\rangle|}{(r-s)^{\al}}\Big(\frac{M^{(1)}_{\gamma}+\al
M^{(2)}_{\gamma}K^*_{s,r}(r-s)^{\beta}}{\Gamma(1-\al)}\Big) \nonumber\\
&     & \mbox{}+\frac{\al M^{(1)}_{\gamma}}{\Gamma(1-\al)} \int_s^r
\frac{|\langle v_r,e_j\rangle-\langle
v_{\theta},e_j\rangle|}{(r-\theta)^{1+\al}} d\theta\Big]\nonumber\\
& = & \sum_{j=1}^n\Big[a_{\gamma,s,r,1}\frac{|\langle
v_r,e_j\rangle|}{(r-s)^{\al}} + \frac{\al
M^{(1)}_{\gamma}}{\Gamma(1-\al)}\int_s^r\frac{|\langle
v_r,e_j\rangle-\langle v_{\theta},e_j\rangle|}{(r-\theta)^{1+\al}}
d\theta\Big]\nonumber\\
& \le & \sum_{j=1}^n\Big[a_{\gamma,s,r,1}|\langle
v_r,e_j\rangle|(r-s)^{-\al} +b_{\gamma,1}\|\langle
v_{\cdot},e_j\rangle\|_{s,t,\beta}(r-s)^{\beta-\al}\Big],
\end{eqnarray*}
where
\begin{equation}
\label{equation.fBM.constt.a(gamma.1)} a_{\gamma,s,r,1} =
\frac{M^{(1)}_{\gamma} +\al
M^{(2)}_{\gamma}K^*_{s,r}(r-s)^{\beta}}{\Gamma(1-\al)},
\end{equation}
and
\begin{equation}
\label{equation.fBm.constt.b(gamma.1)} b_{\gamma,1} = \frac{\al
M^{(1)}_{\gamma}}{(\beta -\al)\Gamma(1-\al)}.
\end{equation}

Note that $a_{\gamma,s,r,1}\le a_{\gamma,s,t,1}$, for $s\le r\le t$.

Writing $a_{s,r,1} = \sum_{\gamma\in I}
a_{\gamma,s,r,1}\|B^H_{\gamma}\|_{\beta,T}$ and $b_1 =
\sum_{\gamma\in I} b_{\gamma,1}\|B^H_{\gamma}\|_{\beta,T}$, and
using the above estimates for the integrands, together with
\eqref{equation.fBM.est.1} and Remark \ref{remark.fBM.alpha.beta},
we have
%and repeating the argument used to obtain \eqref{equation.fBM.flow.est.2} in
%Lemma \ref{lemma.fBM.flow.est}, we get
\begin{eqnarray*}
|\langle (v_t-v_s),e_i\rangle| & \le & k_1(\al,\beta)\int_s^t
\sum_{j=1}^n \Big(a_{s,r,1} |\langle v_r,e_j\rangle|
(r-s)^{-\al}(t-r)^{\al+\beta-1} \\
&     & \mbox{}+b_1\|\langle
v_{\cdot},e_j\rangle\|_{s,t,\beta}(r-s)^{\beta-\al}(t-r)^{\al+\beta-1}
\Big)dr\\
& \le & k_1(\al,\beta) (t-s)^{\al+\beta-1}\int_s^t
\sum_{j=1}^n \Big(a_{s,r,1}|\langle v_r,e_j\rangle| (r-s)^{-\al} \\
&     & \mbox{}+ b_1\|\langle v_{\cdot},e_j\rangle\|_{s,t,\beta}
(r-s)^{\beta-\al}
\Big)dr \\
& \le & k_1(\al,\beta) (t-s)^{\al+\beta-1}\int_s^t \sum_{j=1}^n
\Big(a_{s,r,1}\|\langle v_{\cdot},e_j\rangle\|_{s,t,\iy}
(r-s)^{-\al} \\
&     & \mbox{}+ b_1\|\langle v_{\cdot},e_j\rangle\|_{s,t,\beta}
(r-s)^{\beta-\al}
\Big)dr \\
& \le & k_1(\al,\beta) (t-s)^{\al+\beta-1}\sum_{j=1}^n
\Big(a_{s,t,1}\|\langle v_{\cdot},e_j\rangle\|_{s,t,\iy}
\frac{(t-s)^{1-\al}}{1-\al}\\
&    & \mbox{}+ b_1\|\langle v_{\cdot},e_j\rangle\|_{s,t,\beta}
\frac{(t-s)^{1+\beta-\al}}{1+\beta-\al} \Big)\\
& = & k_1(\al,\beta) \sum_{j=1}^n \Big(a_{s,t,2}\|\langle
v_{\cdot},e_j\rangle\|_{s,t,\iy}
(t-s)^{\beta}\\
&    & \mbox{}+ b_2\|\langle v_{\cdot},e_j\rangle\|_{s,t,\beta}
(t-s)^{2\beta} \Big),
\end{eqnarray*}
where $a_{s,t,2}=a_{s,t,1}(1-\al)^{-1}$ and
$b_2=b_1(1-\al+\beta)^{-1}$.

%Since for any $s\le r\le \theta\le t$, such that $(t-s)$ satisfies
%\eqref{equation.fBM.interval.size.1}, $(\theta-r)$ also satisfies the same, hence
%the above estimate holds for any such $r$ and $\theta$.
Therefore,
\begin{eqnarray*}
\|\langle v_{\cdot},e_i\rangle\|_{s,t,\beta} & = & \sup_{s\le r\le
\theta \le t}\frac{|\langle
(v_{\theta}-v_r),e_i\rangle|}{(\theta-r)^{\beta}}\\
& \le & k_1(\al,\beta) \sum_{j=1}^n \sup_{s\le r\le \theta \le
t}\Big(a_{r,\theta,2}\|\langle v_{\cdot},e_j\rangle\|_{r,\theta,\iy}\\
&    & \mbox{}+ b_2\|\langle v_{\cdot},e_j\rangle\|_{r,\theta,\beta}
(\theta-r)^{\beta} \Big)\\
& \le & k_1(\al,\beta) \sum_{j=1}^n \Big(a_{s,t,2}\|\langle
v_{\cdot},e_j\rangle\|_{s,t,\iy}\\
&    & \mbox{}+ b_2\|\langle v_{\cdot},e_j\rangle\|_{s,t,\beta}
(t-s)^{\beta} \Big).
\end{eqnarray*}

As a consequence of the above estimate we have
\begin{eqnarray}\label{equation.fBM.pre.estimate}
\sum_{i=1}^n\|\langle v_{\cdot},e_i\rangle\|_{s,t,\beta} & \le &
n\;k_1(\al,\beta) \sum_{j=1}^n \Big(a_{s,t,2}\|\langle
v_{\cdot},e_j\rangle\|_{s,t,\iy}\nonumber\\
&    & \mbox{}+ b_2\|\langle v_{\cdot},e_j\rangle\|_{s,t,\beta}
(t-s)^{\beta} \Big).
\end{eqnarray}

For further analysis we shall require that
\begin{eqnarray}\label{equation.fBM.interval.size.2}
(t-s)^{-\beta}> nk_1(\al,\beta)b_2.
\end{eqnarray}
%In view of condition \eqref{equation.fBM.interval.size.1}, the above restraint
%boils down to
%$$c\cdot c_{\al}\tilde{M}^{(1)}_{\al} > b_2.$$
%Equivalently,
%\begin{eqnarray}\label{equation.fBM.interval.size.2}
%c & > & \frac{2(1-2\al)(1-\al)\sum_{\gamma\in
%I}M^{(2)}_{\gamma}\|B^H_{\gamma}\|_{\beta,T}}{(\beta-\al)(1-\al+\beta)\sum_{\gamma\in
%I}M^{(1)}_{\gamma}\|B^H_{\gamma}\|_{\beta,T}}.
%\end{eqnarray}

Thereby, for $(t-s)$ satisfying conditions
\eqref{equation.fBM.interval.size.1} and
\eqref{equation.fBM.interval.size.2}, we can rewrite
\eqref{equation.fBM.pre.estimate} as
$$\sum_{i=1}^n\|\langle v_{\cdot},e_i\rangle\|_{s,t,\beta}\le
n\,k_1(\al,\beta)\,a_{s,t,2}\,\sum_{i=1}^n\frac{\|\langle
v_{\cdot},e_i\rangle\|_{s,t,\iy}}{(1-n\,k_1(\al,\beta)\,b_2\,(t-s)^{\beta})}.$$

Hence,
\begin{eqnarray*}
\sum_{i=1}^n|\langle v_t,e_i\rangle|& \le &
\sum_{i=1}^n\Big(|\langle v_s,e_i\rangle|+
|\langle v_t,e_i\rangle-\langle v_s,e_i\rangle|\Big)\\
& \le & \sum_{i=1}^n\Big(|\langle v_s,e_i\rangle|+
\|\langle v_{\cdot},e_i\rangle\|_{s,t,\beta}(t-s)^{\beta}\Big)\\
& \le & \sum_{i=1}^n\Big(|\langle v_s,e_i\rangle|+
n\,k_1(\al,\beta)\,a_{s,t,2} \frac{\|\langle
v_{\cdot},e_i\rangle\|_{s,t,\iy}(t-s)^{\beta}}{(1
-n\,k_1(\al,\beta)\,b_2\,(t-s)^{\beta})}\Big)
\end{eqnarray*}

Clearly, for any $r\in [s,t]$ we have
$$\sum_{i=1}^n|\langle v_r,e_i\rangle| \le \sum_{i=1}^n
\Big(|\langle v_s,e_i\rangle|+ n\,k_1(\al,\beta)\,a_{s,r,2}
\frac{\|\langle v_{\cdot},e_i\rangle\|_{s,r,\iy}(r-s)^{\beta}}{(1
-n\,k_1(\al,\beta)\,b_2\,(r-s)^{\beta})}\Big).$$

Now using the fact that $s<r<t$, so that $\|\langle
v_{\cdot},e_i\rangle\|_{s,r,\iy}\le \|\langle
v_{\cdot},e_i\rangle\|_{s,t,\iy}$ and $a_{s,r,2}\le a_{s,t,2}$, we
have
\begin{eqnarray}\label{equation.fBM.pre.estimate.1}
\sum_{i=1}^n\|\langle v_{\cdot},e_i\rangle\|_{s,t,\iy} & \le &
\sum_{i=1}^n\Big(|\langle v_s,e_i\rangle|\nonumber\\
&    &\mbox{} + n\,k_1(\al,\beta)\,a_{s,t,2} \frac{\|\langle
v_{\cdot},e_i\rangle\|_{s,t,\iy}(t-s)^{\beta}}{(1
-n\,k_1(\al,\beta)\,b_2\,(t-s)^{\beta})}\Big).
\end{eqnarray}

Finally, we shall require $(t-s)$ to satisfy
\begin{eqnarray}\label{equation.fBM.interval.size.pre.2}
(t-s)^{-\beta} > n\,k_1(\al,\beta)\,[a_{s,t,2}+b_2],
\end{eqnarray}
to allow us to rewrite \eqref{equation.fBM.pre.estimate.1} as
$$\sum_{i=1}^n\|\langle v_{\cdot},e_i\rangle\|_{s,t,\iy}
\Big[1-\frac{n\,k_1(\al,\beta)\,a_{s,t,2}\,(t-s)^{\beta}}{ (1
-n\,k_1(\al,\beta)\,b_2\,(t-s)^{\beta})}\Big] \le
\sum_{i=1}^n|\langle v_s,e_i\rangle|.$$

We shall note that for $(t-s)$ sufficiently small, the inequality
\eqref{equation.fBM.interval.size.pre.2} does hold true, as
$a_{s,t,2}$ is a decreasing function of $(t-s)$.

This, in turn implies,
\begin{eqnarray}
\label{equation.fBM.multiplicative} \sum_{i=1}^n\sup_{0\le r\le
t}|\langle v_r,e_i\rangle| & = & \sum_{i=1}^n\max\{\sup_{0\le r\le
s}|\langle v_r,e_i\rangle| ,
\|\langle v_{\cdot},e_i\rangle\|_{s,t,\iy}\}\nonumber\\
& \le & \sum_{i=1}^n\max \{\sup_{0\le r\le s}|\langle
v_r,e_i\rangle|, \frac{|\langle v_s,e_i\rangle|}{\Big[1-
\frac{nk_1(\al,\beta)a_{s,t,2}(t-s)^{\beta}}{(1
-nk_1(\al,\beta)b_2(t-s)^{\beta})}\Big]}
\}\nonumber\\
& \le & \sum_{i=1}^n\max \{\sup_{0\le r\le s}|\langle
v_r,e_i\rangle|, \frac{\sup_{0\le r\le s}|\langle
v_r,e_i\rangle|}{\Big[1- \frac{nk_1(\al,\beta)a_{s,t,2}
(t-s)^{\beta}}{(1-nk_1(\al,\beta)b_2
(t-s)^{\beta})}\Big]}\}\nonumber\\
& = & \sum_{i=1}^n \frac{\sup_{0\le r\le s}|\langle
v_r,e_i\rangle|}{\Big[1-\frac{nk_1(\al,\beta)a_{s,t,2}
(t-s)^{\beta}}{(1-nk_1(\al,\beta)b_2
(t-s)^{\beta})}\Big]}\nonumber\\
& = & S\sum_{i=1}^n \sup_{0\le r\le s}|\langle v_r,e_i\rangle|,
\end{eqnarray}
where $S=\Big[1- \frac{nk_1(\al\beta)a_{s,t,2}
(t-s)^{\beta}}{(1-nk_1(\al,\beta)b_2(t-s)^{\beta})}\Big]^{-1}.$

Next we divide the interval $[0,T]$ into $p$ pieces of size $\Delta
=(t-s)$, with $\Delta$ being small enough, so that none of the above
estimates fail, and  write $a_{\Delta,2}$ for
$a_{s,t,2}$, as $a_{s,t,2}$ depends on $s,t$ only through the
difference $(t-s)=\Delta$.

More precisely, in view of \eqref{equation.fBM.interval.size.1},
\eqref{equation.fBM.interval.size.2} and
\eqref{equation.fBM.interval.size.pre.2}, we require $\Delta$ to
satisfy
\begin{eqnarray*}
{\Delta}^{-\beta} & > & n\,k_1(\al,\beta)\cdot
\max[c_{\al}\tilde{M}^{(1)}_{\al},\, b_2,\,(a_{\Delta,2}+b_2)]\\
& = & n\,k_1(\al,\beta)\cdot
\max[c_{\al}\tilde{M}^{(1)}_{\al},\,(a_{\Delta,2}+b_2)].
\end{eqnarray*}
For example, we can choose
\begin{equation}
\label{equation.fBM.interval.size.final} \Delta^{-\beta} =
3\,n\,k_1(\al,\beta)\cdot
\max[c_{\al}\tilde{M}^{(1)}_{\al},\,(a_{\Delta,2}+b_2)],
\end{equation}
and thus, for this specific choice of $\Delta$, we have $S\le 2.$

To ensure the existence of such a $\Delta$, we start with
$$\Delta^{-\beta}_0 = 3\,n\,k_1(\al,\beta)\,c_{\al}
\,\tilde{M}^{(1)}_{\al}.$$
Then, if
\begin{equation}\label{equation.fBM.soltn.interval.size}
\Delta^{-\beta}_0 \ge 3\,n\,k_1(\al,\beta)\,(a_{\Delta_0,2}+b_2),
\end{equation}
we choose $\Delta = \Delta_0$. Otherwise we solve the
equation $$\Delta^{-\beta} =
3\,n\,k_1(\al,\beta)\,(a_{\Delta,2}+b_2),$$ in the range $\Delta
\le\Delta_0$. It is easy to see that the solution to this
equation is ensured since the left side increases to infinity as
$\Delta\rightarrow 0$, whereas the right side, which is larger than
the left side at $\Delta=\Delta_0$, decreases as $\Delta$ decreases
to zero.

Using the above notation, and  repeatedly applying the technique
used in \eqref{equation.fBM.multiplicative}, we can write
\begin{eqnarray*}
\sup_{t\in[0,T]}\|v_t\|_1 & = &
\sup_{t\in[0,\;p\Delta]}[\sum_{i=1}^n|\langle v_t,e_i\rangle|]\\
& \le & \sum_{i=1}^n\sup_{t\in[0,\;p\Delta]}|\langle v_t,e_i\rangle|\\
& \le & S^p\sum_{i=1}^n |\langle v,e_i\rangle|,
\end{eqnarray*}
where
\begin{eqnarray*}
p \ = \ \frac{T}{\Delta}
& = & T\Big(3nk_1(\al,\beta)\cdot
\max[c_{\al}\tilde{M}^{(1)}_{\al},\,(a_{\Delta,2}+b_2)]\Big)^{1/\beta}\\
& = & T\;C_T,
\end{eqnarray*}
and
$$C_T = \Big(3n\,k_1(\al,\beta)\cdot
\max[(c_{\al}\tilde{M}^{(1)}_{\al}),\,(a_{\Delta,2}+b_2)]
\Big)^{1/\beta}.$$

Since we have all the appropriate notation at hand,
we now take a moment off the proof of Theorem {\ref{theorem.fBM.growth.vector}
to complete the remaining issues in the proof of
Lemma \ref{lemma.fBM.flow.est}

\textbf{Proof of Lemma {\ref{lemma.fBM.flow.est}} (continued):}
To prove the final claim of Lemma \ref{lemma.fBM.flow.est}, note
that for specific choice $(t-s)=\Delta$, together with
\eqref{equation.fBM.constt.M(alpha)} and
\eqref{equation.fBM.flow.pre.est} we have
\begin{eqnarray*}
K^*_{s,t} & = & \frac{K_{s,t}}{(1-2\al)}\\
& \le & \frac{3}{2(1-2\al)}\cdot c_{\al}\,n\,k_1(\al,\beta)M_{\al}\\
& = & \frac{3}{2(1-2\al)}\cdot c_{\al}\,n\,k_1(\al,\beta)
\sum_{\gamma\in I}
\frac{M_{\gamma}\|B^H_{\gamma}\|_{\beta,T}}{1-\al},
\end{eqnarray*}
and so there exists a constant $K(\al,\beta)$, dependent only on
$\al$ and $\beta$, such that
\begin{eqnarray*}
K^*_{s,t}(t-s)^{\beta} & \le & K(\al,\beta)\frac{\sum_{\gamma\in I}
M_{\gamma}\|B^H_{\gamma}\|_{\beta,T}}{2\sum_{\gamma\in
I}M^{(1)}_{\gamma}\|B^H_{\gamma}\|_{\beta,T}}\\
& \le & K(\al,\beta)\sum_{\gamma\in
I}\frac{M_{\gamma}}{M^{(1)}_{\gamma}}.
\end{eqnarray*}
Consequently, $a_{\Delta,2}$ can also be bounded above by a constant
$a_2$, hence we shall replace $a_{\Delta,2}$ by $a_2$, in the
following discussion. \qed \\

Returning to  the proof of Theorem {\ref{theorem.fBM.growth.vector},
note that
$$(C_T)^{\beta} \le 3n\,c_{\al}\,k_1(\al,\beta) \sum_{\gamma\in I}
(\tilde{M}^{(1)}_{\al,\gamma}+a_{2,\gamma}+b_{2,\gamma})
\|B^H_{\gamma}\|_{\beta,T},$$ where $\tilde{M}^{(1)}_{\al,\gamma}$,
$a_{2,\gamma}$, and $b_{2,\gamma}$ are the coefficients of
$\|B^H_{\gamma}\|_{\beta,T}$ in the constants
$\tilde{M}^{(1)}_{\al}$, $a_2$ and $b_2$, respectively.

Now using the bound on $S$ available due to the specific choice of
$\Delta$ completes the proof. \qed \\

The estimates in Theorem \ref{theorem.fBM.growth.vector} in turn
imply similar bounds on the Hausdorff measure of the $m$-dimensional
manifold $M_t$, evolving under the flow $\Phi_t$. More precisely,
let $\{v_i^x\}_{i=1}^m$ be an orthonormal basis of the tangent space
$T_xM$, at the point $x\in M$. Then, writing, as usual, $\calH_m(M_t)$ for
$m$-th Hausdorff measure of $M_t$, we have the following
result.

\begin{theorem}
\label{theorem.fBM.growth.Hausdorff.measure} Let $M$ be a $C^2$, 
$m$-dimensional manifold, evolving under the flow $\Phi_t$ defined in
\eqref{equation.fBM.flow.candidate}. Then under the conditions
$(A1)-(A4)$, and for $\al= 1-H+\delta$, $\beta=H-\eps$, such that
$(1-H)<\al<1/2$ and $\delta>\eps$, there exists a constant $c_1$,
and a random variable $C_{1,T}$, such that
$$\sup_{t\in[0,T]}\calH_m(M_t)\;\le\; c_1\,\calH_m(M) \;2^{C_{1,T}\,T}.$$
Here  $C_{1,T}$ depends on $\al$, $\beta$, $n$, $I$,
and $\{\|B^H_{\gamma}\|_{\beta,T}, M_{\gamma}, M^{(1)}_{\gamma},
M^{(2)}_{\gamma}\}_{\gamma\in I}$,
and satsifies
 $$E[C_{1,T}]^{\beta} \le C_1\cdot E [\|B^H\|_{\beta,T}],$$
with the constant $C_1$ dependent only on $\al$, $\beta$, $n$, $|I|$ and
$\{M_{\gamma}, M^{(1)}_{\gamma}, M^{(2)}_{\gamma}\}_{\gamma\in I}$.
\end{theorem}

\textbf{Proof:} Consider the pushforwards $\{v_{i,t}^x\}_{i=1}^m$ of the
tangent vectors $\{v_i^x\}_{i=1}^m$ under the flow $\Phi_t$. Then by
using a simple formula for change of variables on a manifold, we have
\begin{eqnarray*}
\calH_m(M_t) & = & \int_{M_t} \calH_m(dy)\\
& = & \int_{M} \|\al^x(t)\| \calH_m(dx),
\end{eqnarray*}
where $\|\al^x(t)\|=\sqrt{|\det(\langle
v_{i,t}^x,v_{j,t}^x\rangle)|}$. By the Cauchy-Schwartz inequality we
know that
$$\langle v_{i,t}^x,v_{j,t}^x\rangle \le \|v_{i,t}^x\|_2\;\|v_{i,t}^x\|_2.$$
Therefore, using Theorem \ref{theorem.fBM.growth.vector} and the
above expression, we obtain
\begin{eqnarray*}
\sup_{t\in[0,T]}\|\al^x(t)\| & \le &
m!(\sup_{t\in[0,T]}\|v_{i,t}^x\|)^m\\
& \le & c\,m!\,2^{m\,T\,C_T},
\end{eqnarray*}
which proves the required result. \qed

Finally, note that for the case corresponding to $|I|=1$, the above
simplifies to
$$\sup_{t\in[0,T]}\calH_m(M_t)\;\le\; c_1\;2^{C_1\;T\;\|B^H\|^{1/\beta}_{\beta,T}},$$
for some constants $c_1$ and $C_1$ dependent only on the various
Lipschitz coefficients of the vector field and its partial
derivatives, and this proves Theorem
\ref{theorem.growth.Hausdorff.1}.

{\smc Sreekar Vadlamani:} Faculty of Industrial Engineering and
Management, Technion - Israel Institute of Technology, Haifa, Israel
\newline sreekar.vadlamani@gmail.com

\end{document}